# DISCUSSION OF "LEAST ANGLE REGRESSION" BY EFRON ET AL.

By David Madigan and Greg Ridgeway

*Rutgers University and Avaya Labs Research, and RAND*

Algorithms for simultaneous shrinkage and selection in regression and classification provide attractive solutions to knotty old statistical challenges. Nevertheless, as far as we can tell, Tibshirani's Lasso algorithm has had little impact on statistical practice. Two particular reasons for this may be the relative inefficiency of the original Lasso algorithm and the relative complexity of more recent Lasso algorithms [e.g., Osborne, Presnell and Turlach (2000)]. Efron, Hastie, Johnstone and Tibshirani have provided an efficient, simple algorithm for the Lasso as well as algorithms for stagewise regression and the new least angle regression. As such this paper is an important contribution to statistical computing.

**1. Predictive performance.** The authors say little about predictive performance issues. In our work, however, the relative out-of-sample predictive performance of LARS, Lasso and Forward Stagewise (and variants thereof) takes center stage. Interesting connections exist between boosting and stagewise algorithms so predictive comparisons with boosting are also of interest.

The authors present a simple $C_p$ statistic for LARS. In practice, a cross-validation (CV) type approach for selecting the degree of shrinkage, while computationally more expensive, may lead to better predictions. We considered this using the LARS software. Here we report results for the authors' diabetes data, the Boston housing data and the Servo data from the UCI Machine Learning Repository. Specifically, we held out 10% of the data and chose the shrinkage level using either $C_p$ or nine-fold CV using 90% of the data. Then we estimated mean square error (MSE) on the 10% hold-out sample. Table 1 shows the results for main-effects models.

Table 1 exhibits two particular characteristics. First, as expected, Stagewise, LARS and Lasso perform similarly. Second, $C_p$ performs as well as cross-validation; if this holds up more generally, larger-scale applications will want to use $C_p$ to select the degree of shrinkage.

TABLE 1
*Stagewise, LARS and Lasso mean square error predictive performance, comparing cross-validation with $C_p$*

|  | **Diabetes** | | | **Boston** | | | **Servo** | |
|---|---|---|---|---|---|---|---|---|
|  | CV | $C_p$ |  | CV | $C_p$ |  | CV | $C_p$ |
| Stagewise | 3083 | 3082 | Stagewise | 25.7 | 25.8 | Stagewise | 1.33 | 1.32 |
| LARS | 3080 | 3083 | LARS | 25.5 | 25.4 | LARS | 1.33 | 1.30 |
| Lasso | 3083 | 3082 | Lasso | 25.8 | 25.7 | Lasso | 1.34 | 1.31 |

Table 2 presents a reanalysis of the same three datasets but now considering five different models: least squares; LARS using cross-validation to select the coefficients; LARS using $C_p$ to select the coefficients and allowing for two-way interactions; least squares boosting fitting only main effects; least squares boosting allowing for two-way interactions. Again we used the authors' LARS software and, for the boosting results, the gbm package in R [Ridgeway (2003)]. We evaluated all the models using the same cross-validation group assignments.

A plain linear model provides the best out-of-sample predictive performance for the diabetes dataset. By contrast, the Boston housing and Servo data exhibit more complex structure and models incorporating higher-order structure do a better job.

While no general conclusions can emerge from such a limited analysis, LARS seems to be competitive with these particular alternatives. We note, however, that for the Boston housing and Servo datasets Breiman (2001) reports substantially better predictive performance using random forests.

TABLE 2
*Predictive performance of competing methods: LM is a main-effects linear model with least squares fitting; LARS is least angle regression with main effects and CV shrinkage selection; LARS two-way $C_p$ is least angle regression with main effects and all two-way interactions, shrinkage selection via $C_p$; GBM additive and GBM two-way use least squares boosting, the former using main effects only, the latter using main effects and all two-way interactions; MSE is mean square error on a 10% holdout sample; MAD is mean absolute deviation*

|  | **Diabetes** | | **Boston** | | **Servo** | |
|---|---|---|---|---|---|---|
|  | MSE | MAD | MSE | MAD | MSE | MAD |
| LM | 3000 | 44.2 | 23.8 | 3.40 | 1.28 | 0.91 |
| LARS | 3087 | 45.4 | 24.7 | 3.53 | 1.33 | 0.95 |
| LARS two-way $C_p$ | 3090 | 45.1 | 14.2 | 2.58 | 0.93 | 0.60 |
| GBM additive | 3198 | 46.7 | 16.5 | 2.75 | 0.90 | 0.65 |
| GBM two-way | 3185 | 46.8 | 14.1 | 2.52 | 0.80 | 0.60 |

**2. Extensions to generalized linear models.** The minimal computational complexity of LARS derives largely from the squared error loss function. Applying LARS-type strategies to models with nonlinear loss functions will require some form of approximation. Here we consider LARS-type algorithms for logistic regression.

Consider the logistic log-likelihood for a regression function $f(\mathbf{x})$ which will be linear in $\mathbf{x}$:

$$(2.1) \quad \ell(f) = \sum_{i=1}^{N} y_i f(\mathbf{x}_i) - \log(1 + \exp(f(\mathbf{x}_i))).$$

We can initialize $f(\mathbf{x}) = \log(\bar{y}/(1-\bar{y}))$. For some $\alpha$ we wish to find the covariate $x_j$ that offers the greatest improvement in the logistic log-likelihood, $\ell(f(\mathbf{x}) + \mathbf{x}_j^t \alpha)$. To find this $\mathbf{x}_j$ we can compute the directional derivative for each $j$ and choose the maximum,

$$(2.2) \quad j^* = \arg\max_j \left| \frac{d}{d\alpha} \ell(f(\mathbf{x}) + \mathbf{x}_j^t \alpha) \right|_{\alpha=0}$$

$$(2.3) \quad = \arg\max_j \left| \mathbf{x}_j^t \left( \mathbf{y} - \frac{1}{1 + \exp(-f(\mathbf{x}))} \right) \right|.$$

Note that as with LARS this is the covariate that is most highly correlated with the residuals. The selected covariate is the first member of the active set, $A$. For $\alpha$ small enough (2.3) implies

$$(2.4) \quad (s_{j^*} \mathbf{x}_{j^*} - s_j \mathbf{x}_j)^t \left( y - \frac{1}{1 + \exp(-f(\mathbf{x}) - \mathbf{x}_{j^*}^t \alpha)} \right) \geq 0$$

for all $j \in A^C$, where $s_j$ indicates the sign of the correlation as in the LARS development. Choosing $\alpha$ to have the largest magnitude while maintaining the constraint in (2.4) involves a nonlinear optimization. However, linearizing (2.4) yields a fairly simple approximate solution. If $x_2$ is the variable with the second largest correlation with the residual, then

$$(2.5) \quad \hat{\alpha} = \frac{(s_{j^*} \mathbf{x}_{j^*} - s_2 \mathbf{x}_2)^t (y - p(\mathbf{x}))}{(s_{j^*} \mathbf{x}_{j^*} - s_2 \mathbf{x}_2)^t (p(\mathbf{x})(1 - p(\mathbf{x})) \mathbf{x}_{j^*})}.$$

The algorithm may need to iterate (2.5) to obtain the exact optimal $\hat{\alpha}$. Similar logic yields an algorithm for the full solution.

We simulated $N = 1000$ observations with 10 independent normal covariates $\mathbf{x}_i \sim N_{10}(\mathbf{0}, \mathbf{I})$ with outcomes $Y_i \sim \text{Bern}(1/(1 + \exp(-\mathbf{x}_i^t \beta)))$, where $\beta \sim N_{10}(0, 1)$. Figure 1 shows a comparison of the coefficient estimates using Forward Stagewise and the Least Angle method of estimating coefficients, the final estimates arriving at the MLE. While the paper presents LARS for squared error problems, the Least Angle approach seems applicable to a wider family of models. However, an out-of-sample evaluation of predictive performance is essential to assess the utility of such a modeling strategy.

Specifically for the Lasso, one alternative strategy for logistic regression is to use a quadratic approximation for the log-likelihood. Consider the Bayesian version of Lasso with hyperparameter $\gamma$ (i.e., the penalized rather than constrained version of Lasso):

$$\log f(\boldsymbol{\beta}|y_1,\ldots,y_n)$$
$$\propto \sum_{i=1}^{n}\log(y_i\Lambda(\mathbf{x}_i\boldsymbol{\beta})+(1-y_i)(1-\Lambda(\mathbf{x}_i\boldsymbol{\beta})))+d\log\left(\frac{\gamma^{1/2}}{2}\right)-\gamma^{1/2}\sum_{i=1}^{d}|\beta_i|$$
$$\approx \left(\sum_{i=1}^{n}a_i(\mathbf{x}_i\boldsymbol{\beta})^2+b_i(\mathbf{x}_i\boldsymbol{\beta})+c_i\right)+d\log\left(\frac{\gamma^{1/2}}{2}\right)-\gamma^{1/2}\sum_{i=1}^{d}|\beta_i|,$$

where $\Lambda$ denotes the logistic link, $d$ is the dimension of $\boldsymbol{\beta}$ and $a_i, b_i$ and $c_i$ are Taylor coefficients. Fu's elegant coordinatewise "Shooting algorithm" [Fu (1998)], can optimize this target starting from either the least squares solution or from zero. In our experience the shooting algorithm converges rapidly.

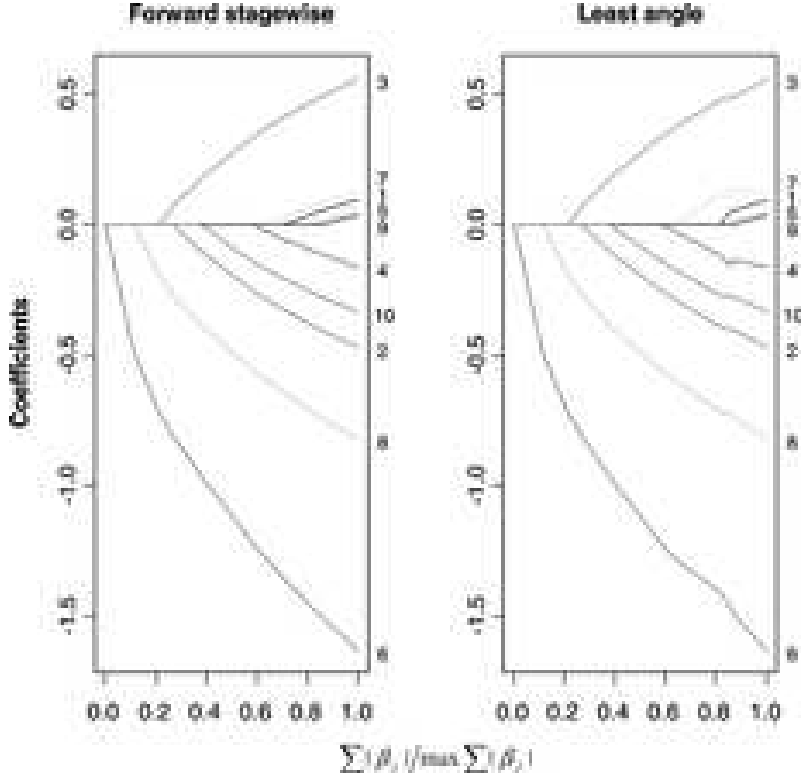

FIG. 1. *Comparison of coefficient estimation for Forward Stagewise and Least Angle Logistic Regression.*

## REFERENCES

BREIMAN, L. (2001). Random forests. Available at ftp:// ftp.stat.berkeley.edu/pub/users/ breiman/randomforest2001.pdf.

FU, W. J. (1998). Penalized regressions: The Bridge versus the Lasso. *J. Comput. Graph. Statist.* **7** 397–416.

OSBORNE, M. R., PRESNELL, B. and TURLACH, B. A. (2000). A new approach to variable selection in least squares problems. *IMA J. Numer. Anal.* **20** 389–403. MR1646710

RIDGEWAY, G. (2003). GBM 0.7-2 package manual. Available at http:// cran.r-project. org/doc/packages/gbm.pdf.

AVAYA LABS RESEARCH
AND
DEPARTMENT OF STATISTICS
RUTGERS UNIVERSITY
PISCATAWAY, NEW JERSEY 08855
USA
E-MAIL: madigan@stat.rutgers.edu

RAND STATISTICS GROUP
SANTA MONICA, CALIFORNIA 90407-2138
USA
E-MAIL: gregr@rand.org